\documentclass[12pt,twoside,reqno]{amsart}
\usepackage{pifont}
\usepackage{amsmath}
\usepackage{amsthm}
\usepackage{amsfonts}
\usepackage{amssymb}
\usepackage{latexsym}
\usepackage{amstext}
\usepackage{array}
\usepackage{graphicx}

\date{}
\allowdisplaybreaks[4] \footskip=15pt
\renewcommand{\uppercasenonmath}[1]{}

\newtheorem{thm}[subsection]{Theorem}
\newtheorem{cor}[subsection]{Corollary }
\newtheorem{Def}[subsection]{Definition}
\newtheorem{lem}[subsection]{Lemma}
\newtheorem{remark}[subsection]{Remark}

\newtheorem{prop}[subsection]{Proposition}
\newtheorem{exm}[subsection]{Example}
  
\newcommand{\bthm}{\begin{thm} }
\newcommand{\ethm}{\end{thm} }
\newcommand{\bpro}{\begin{prop}}
\newcommand{\epro}{\end{prop}}
\newcommand{\bdf}{\begin{Def}}
\newcommand{\edf}{\end{Def}}
\newcommand{\bexm}{\begin{exm}}
\newcommand{\eexm}{\end{exm}}
\newcommand{\blem}{\begin{lem}}
\newcommand{\elem}{\end{lem}}
\newcommand{\bpf}{\begin{proof}}
\newcommand{\epf}{\end{proof}}
\newcommand{\bcor}{\begin{cor}}
\newcommand{\ecor}{\end{cor}}
\newcommand{\ba}{\begin{array}}
\newcommand{\ea}{\end{array}}
\newcommand{\bea}{\begin{eqnarray}}
\newcommand{\eea}{\end{eqnarray}}

\newcommand{\brem}{\begin{remark}}
\newcommand{\erem}{\end{remark}}

\begin{document}
\begin{center}
{\large  \bf Higher $K$-Groups of Smooth Projective Curves Over
Finite Fields }\\

 \vskip 0.8cm
 {\small  {\ding {172}} Qingzhong Ji \ \   {\ding {173}}} Hourong Qin\\
{\small Department of Mathematics, Nanjing University, Nanjing
210093, P.R.China}\\
{\small E-mail: \ding {172} qingzhji@nju.edu.cn \;\;\; \ding {173}
hrqin@nju.edu.cn}
\end{center}

\newsavebox{\tablebox}
\begin{lrbox}{\tablebox}
\begin{tabular}{|c|c|c|c|c|c|l|}
\multicolumn{7}{c}{\bf TABLE\;\;\; I(1)}\\
\multicolumn{7}{c}{(\text{of\; all\; elliptic\; curves}\; E
\;\text{defined\; over}\;
$\mathbb{F}=\mathbb{Z}/3\mathbb{Z}$)}\\

\multicolumn{7}{c}{}\\

\hline
 No. &\multicolumn{1}{|c|}{$E$} & \multicolumn{1}{|c|}{\text{roots}}&
\multicolumn{1}{|c|}{$E(\mathbb{F})$}&
\multicolumn{1}{|c|}{$K_2(E)$}&
\multicolumn{1}{|c|}{$\lambda(l)$}& \multicolumn{1}{|c|}{$l$-{\text{Sylow}}\;\;\text{subgroups}}\\
\hline

$1$& {$ y^2=x^3-x-1$} & $\frac{3\pm\sqrt{-3}}{2}$&
\{$O$\}&$\frac{\mathbb{Z}}{19\mathbb{Z}}$& $\lambda(19)=1$&{
$K_2(19^m)(19)\cong \frac{\mathbb{Z}}{19^{m+1}\mathbb{Z}},\;m\geq 1$}\\
\hline

 {}&{}&{}&{}&{}&
\multicolumn{1}{|c|}{$\lambda(2)=2$}&\multicolumn{1}{|c|}{$K_2(2^m)(2)\cong
\frac{\mathbb{Z}}{2^{m+1}\mathbb{Z}}\times\frac{\mathbb{Z}}{2^{m}\mathbb{Z}},\;m\geq 2$}\\
\cline{6-7} \raisebox{1.6ex}[0pt]{$2$}&
\raisebox{1.6ex}[0pt]{$y^2=x^3-x^2-1$} &
\raisebox{1.6ex}[0pt]{$1\pm\sqrt{-2}$}&
\raisebox{1.6ex}[0pt]{$\frac{\mathbb{Z}}{2\mathbb{Z}}$}&
\raisebox{1.6ex}[0pt]{$\frac{\mathbb{Z}}{22\mathbb{Z}}$}&
\multicolumn{1}{|c|}{$\lambda(11)=1$}&\multicolumn{1}{|c|}{$K_2(11^m)(11)\cong
\mathbb{Z}/11^{m+1}\mathbb{Z},\;m\geq 1$}\\
\hline

$3$&
$y^2=x^3+x^2-1$&$\frac{1\pm\sqrt{-11}}{2}$&$\frac{\mathbb{Z}}{3\mathbb{Z}}$&
$\frac{\mathbb{Z}}{25\mathbb{Z}}$&$\lambda(5)=1$&$K_2(5^m)(5)\cong
\frac{\mathbb{Z}}{5^{m+2}\mathbb{Z}},\;m\geq 0$\\
\hline

{}&{}&{}&{}&{}&
\multicolumn{1}{|c|}{$\lambda(2)=2$}&\multicolumn{1}{|c|}{$K_2(2^{m})(2)\cong
\frac{\mathbb{Z}}{2^{m+1}\mathbb{Z}}\times\frac{\mathbb{Z}}{2^{m+1}\mathbb{Z}},\;m\geq 1$}\\
\cline{6-7} \raisebox{1.6ex}[0pt]{$4$}&
\raisebox{1.6ex}[0pt]{$y^2=x^3+x$} &
\raisebox{1.6ex}[0pt]{$\pm\sqrt{-3}$}&
\raisebox{1.6ex}[0pt]{$\frac{\mathbb{Z}}{4\mathbb{Z}}$}&
\raisebox{1.6ex}[0pt]{$\frac{\mathbb{Z}}{28\mathbb{Z}}$}&
\multicolumn{1}{|c|}{$\lambda(7)=1$}&\multicolumn{1}{|c|}{$K_2(7^m)(7)\cong
\frac{\mathbb{Z}}{7^{m+1}\mathbb{Z}},\;m\geq 0$}\\
\hline

{}&{}&{}&{}&{}&
\multicolumn{1}{|c|}{$\lambda(2)=2$}&\multicolumn{1}{|c|}{$K_2(2^{m})(2)\cong
\frac{\mathbb{Z}}{2^{m+1}\mathbb{Z}}\times\frac{\mathbb{Z}}{2^{m+1}\mathbb{Z}},\;m\geq 0$}\\
\cline{6-7} \raisebox{1.6ex}[0pt]{$5$}&
\raisebox{1.6ex}[0pt]{$y^2=x^3-x$} &
\raisebox{1.6ex}[0pt]{$\pm\sqrt{-3}$}&
\raisebox{1.6ex}[0pt]{$\frac{\mathbb{Z}}{2\mathbb{Z}}\times\frac{\mathbb{Z}}{2\mathbb{Z}}$}&
\raisebox{1.6ex}[0pt]{$\frac{\mathbb{Z}}{2\mathbb{Z}}\times\frac{\mathbb{Z}}{14\mathbb{Z}}$}&
\multicolumn{1}{|c|}{$\lambda(7)=1$}&\multicolumn{1}{|c|}{$K_2(7^m)(7)\cong
\frac{\mathbb{Z}}{7^{m+1}\mathbb{Z}},\;m\geq 0$}\\
\hline

$6$&
$y^2=x^3-x^2+1$&$\frac{-1\pm\sqrt{-11}}{2}$&$\frac{\mathbb{Z}}{5\mathbb{Z}}$&
$\frac{\mathbb{Z}}{31\mathbb{Z}}$&$\lambda(31)=1$&$K_2(31^m)(31)\cong
\frac{\mathbb{Z}}{31^{m+1}\mathbb{Z}},\;m\geq 0$\\
\hline

 {}&{}&{}&{}&{}&
\multicolumn{1}{|c|}{$\lambda(2)=2$}&\multicolumn{1}{|c|}{$K_2(2^m)(2)\cong
\frac{\mathbb{Z}}{2^{m+1}\mathbb{Z}}\times\frac{\mathbb{Z}}{2^{m}\mathbb{Z}},\;m\geq 2$}\\
\cline{6-7} \raisebox{1.6ex}[0pt]{$7$}&
\raisebox{1.6ex}[0pt]{$y^2=x^3+x^2+1$} &
\raisebox{1.6ex}[0pt]{$-1\pm\sqrt{-2}$}&
\raisebox{1.6ex}[0pt]{$\frac{\mathbb{Z}}{6\mathbb{Z}}$}&
\raisebox{1.6ex}[0pt]{$\frac{\mathbb{Z}}{34\mathbb{Z}}$}&
\multicolumn{1}{|c|}{$\lambda(17)=1$}&\multicolumn{1}{|c|}{$K_2(17^m)(17)\cong
\frac{\mathbb{Z}}{17^{m+1}\mathbb{Z}},\;m\geq 0$}\\
\hline

$8$& { $y^2=x^3-x+1$} & $\frac{-3\pm\sqrt{-3}}{2}$&
$\frac{\mathbb{Z}}{7\mathbb{Z}}$&$\frac{\mathbb{Z}}{37\mathbb{Z}}$&
$\lambda(37)=1$&${ K_2(37^m)(37)\cong
\frac{\mathbb{Z}}{37^{m+1}\mathbb{Z}},\;m \geq 0}$\\
\hline

\multicolumn{7}{|c|}{\text{The\; elliptic\; curves\; of\; numbers}\;
$n$\; \text{and}\; $9-n$\; \text{are\; twisted.}}\\\hline

\end{tabular}
\end{lrbox}
\resizebox{5.5in}{!}{\usebox{\tablebox}}

\begin{lrbox}{\tablebox}
\begin{tabular}{|c|c|c|c|c|c|c|c|}
\multicolumn{8}{c}{\bf TABLE\;\;\; I(2)}\\
\multicolumn{8}{c}{(\text{of\; all\; elliptic\; curves}\; E
\;\text{defined\; over}\;
$\mathbb{F}=\mathbb{Z}/3\mathbb{Z}$)}\\

\multicolumn{8}{c}{}\\
\hline

No. &\multicolumn{1}{|c|}{$E$} &\multicolumn{1}{|c|}{$K_2(E)$}&
\multicolumn{1}{|c|}{$K_4(E)$}& \multicolumn{1}{|c|}{$K_6(E)$}&
\multicolumn{1}{|c|}{$K_8(E)$}& \multicolumn{1}{|c|}{$K_{10}(E)$}&
\multicolumn{1}{|c|}{$K_{12}(E)$}\\
\hline

$1$ & $y^2=x^3-x-1$& $\frac{\mathbb{Z}}{19\mathbb{Z}}$&
$\frac{\mathbb{Z}}{217\mathbb{Z}}$&
$\frac{\mathbb{Z}}{2107\mathbb{Z}}$&
$\frac{\mathbb{Z}}{19441\mathbb{Z}}$&
$\frac{\mathbb{Z}}{176419\mathbb{Z}}$&
$\frac{\mathbb{Z}}{1592137\mathbb{Z}}$\\
\hline

$2$ & $y^2=x^3-x^2-1$& $\frac{\mathbb{Z}}{22\mathbb{Z}}$&
$\frac{\mathbb{Z}}{226\mathbb{Z}}$&
$\frac{\mathbb{Z}}{2134\mathbb{Z}}$&
$\frac{\mathbb{Z}}{19552\mathbb{Z}}$&
$\frac{\mathbb{Z}}{176662\mathbb{Z}}$&
$\frac{\mathbb{Z}}{1592866\mathbb{Z}}$\\
\hline

$3$ & $y^2=x^3+x^2-1$& $\frac{\mathbb{Z}}{25\mathbb{Z}}$&
$\frac{\mathbb{Z}}{235\mathbb{Z}}$&
$\frac{\mathbb{Z}}{2161\mathbb{Z}}$&
$\frac{\mathbb{Z}}{19603\mathbb{Z}}$&
$\frac{\mathbb{Z}}{176905\mathbb{Z}}$&
$\frac{\mathbb{Z}}{1593595\mathbb{Z}}$\\
\hline

$4$ & $y^2=x^3+x$& $\frac{\mathbb{Z}}{28\mathbb{Z}}$&
$\frac{\mathbb{Z}}{244\mathbb{Z}}$&
$\frac{\mathbb{Z}}{2188\mathbb{Z}}$&
$\frac{\mathbb{Z}}{19684\mathbb{Z}}$&
$\frac{\mathbb{Z}}{177148\mathbb{Z}}$&
$\frac{\mathbb{Z}}{1594324\mathbb{Z}}$\\
\hline

$5$ & $y^2=x^3-x$&
$\frac{\mathbb{Z}}{2\mathbb{Z}}\times\frac{\mathbb{Z}}{14\mathbb{Z}}$&
$\frac{\mathbb{Z}}{2\mathbb{Z}}\times\frac{\mathbb{Z}}{122\mathbb{Z}}$&
$\frac{\mathbb{Z}}{2\mathbb{Z}}\times\frac{\mathbb{Z}}{1094\mathbb{Z}}$&
$\frac{\mathbb{Z}}{2\mathbb{Z}}\times\frac{\mathbb{Z}}{9842\mathbb{Z}}$&
$\frac{\mathbb{Z}}{2\mathbb{Z}}\times\frac{\mathbb{Z}}{88574\mathbb{Z}}$&
$\frac{\mathbb{Z}}{2\mathbb{Z}}\times\frac{\mathbb{Z}}{797162\mathbb{Z}}$\\
\hline

$6$ & $y^2=x^3-x^2+1$& $\frac{\mathbb{Z}}{31\mathbb{Z}}$&
$\frac{\mathbb{Z}}{253\mathbb{Z}}$&
$\frac{\mathbb{Z}}{2215\mathbb{Z}}$&
$\frac{\mathbb{Z}}{19765\mathbb{Z}}$&
$\frac{\mathbb{Z}}{177391\mathbb{Z}}$&
$\frac{\mathbb{Z}}{1595053\mathbb{Z}}$\\
\hline

$7$ & $y^2=x^3+x^2+1$& $\frac{\mathbb{Z}}{34\mathbb{Z}}$&
$\frac{\mathbb{Z}}{262\mathbb{Z}}$&
$\frac{\mathbb{Z}}{2242\mathbb{Z}}$&
$\frac{\mathbb{Z}}{19846\mathbb{Z}}$&
$\frac{\mathbb{Z}}{177634\mathbb{Z}}$&
$\frac{\mathbb{Z}}{1595782\mathbb{Z}}$\\
\hline

$8$ & $y^2=x^3-x+1$& $\frac{\mathbb{Z}}{37\mathbb{Z}}$&
$\frac{\mathbb{Z}}{271\mathbb{Z}}$&
$\frac{\mathbb{Z}}{2269\mathbb{Z}}$&
$\frac{\mathbb{Z}}{19927\mathbb{Z}}$&
$\frac{\mathbb{Z}}{177877\mathbb{Z}}$&
$\frac{\mathbb{Z}}{1596511\mathbb{Z}}$\\

\hline \multicolumn{8}{|c|}{\text{The\; elliptic\; curves\; of\;
numbers}\; $n$\;\text{ and}\; $9-n$\; \text{ are\;
twisted}.}\\\hline

\end{tabular}
\end{lrbox}
\resizebox{5.5 in}{2.0 in}{\usebox{\tablebox}}

\begin{lrbox}{\tablebox}
\begin{tabular}{|c|c|c|c|c|c|c|}
\multicolumn{7}{c}{\bf TABLE\;\;\; II(1)}\\
\multicolumn{7}{c}{(\text{of\; all\; elliptic\; curves}\; E
\;\text{defined\; over}\;
$\mathbb{F}=\mathbb{Z}/5\mathbb{Z}$)}\\

\multicolumn{7}{c}{}\\
\hline

 No. &\multicolumn{1}{|c|}{$E$} & \multicolumn{1}{|c|}{\text{roots}}&
\multicolumn{1}{|c|}{$E(\mathbb{F})$}&
\multicolumn{1}{|c|}{$K_2(E)$}&
\multicolumn{1}{|c|}{$\lambda(l)$}& \multicolumn{1}{|c|}{ $l$-{\text{Sylow}}\;\;\text{subgroups}}\\
\hline

{}&{}&{}&{}&{}&
\multicolumn{1}{|c|}{$\lambda(2)=2$}&\multicolumn{1}{|c|}{$K_2(2^m)(2)\cong
\frac{\mathbb{Z}}{2^{m+2}\mathbb{Z}}\times\frac{\mathbb{Z}}{2^{m+3}\mathbb{Z}},\;m\geq 2$}\\
\cline{6-7} \raisebox{1.6ex}[0pt]{$1$}&
\raisebox{1.6ex}[0pt]{$y^2=x^3+2x$} &
\raisebox{1.6ex}[0pt]{$2\pm\sqrt{-1}$}&
\raisebox{1.6ex}[0pt]{$\frac{\mathbb{Z}}{2\mathbb{Z}}$}&
\raisebox{1.6ex}[0pt]{$\frac{\mathbb{Z}}{106\mathbb{Z}}$}&
\multicolumn{1}{|c|}{$\lambda(53)=1$}&\multicolumn{1}{|c|}{$K_2(53^m)(53)\cong
\frac{\mathbb{Z}}{53^{m+1}\mathbb{Z}},\;m\geq 0$}\\
\hline

{}&{}&{}&{}&{}&
\multicolumn{1}{|c|}{$\lambda(3)=1$}&\multicolumn{1}{|c|}{$K_2(3^m)(3)\cong
\frac{\mathbb{Z}}{3^{m+1}\mathbb{Z}},\;m\geq 0$}\\
\cline{6-7} \raisebox{1.6ex}[0pt]{$2$}&
\raisebox{1.6ex}[0pt]{$y^2=x^3-x+2$} &
\raisebox{1.6ex}[0pt]{$\frac{3\pm\sqrt{-11}}{2}$}&
\raisebox{1.6ex}[0pt]{$\frac{\mathbb{Z}}{3\mathbb{Z}}$}&
\raisebox{1.6ex}[0pt]{$\frac{\mathbb{Z}}{111\mathbb{Z}}$}&
\multicolumn{1}{|c|}{$\lambda(37)=1$}&\multicolumn{1}{|c|}{$K_2(37^m)(37)\cong
\frac{\mathbb{Z}}{37^{m+1}\mathbb{Z}},\;m\geq 0$}\\
\hline

{}&{}&{}&{}&{}&
\multicolumn{1}{|c|}{$\lambda(2)=2$}&\multicolumn{1}{|c|}{$K_2(2^m)(2)\cong
\frac{\mathbb{Z}}{2^{m+1}\mathbb{Z}}\times\frac{\mathbb{Z}}{2^{m+2}\mathbb{Z}},\;m\geq 1$}\\
\cline{6-7} \raisebox{1.6ex}[0pt]{$3$}&
\raisebox{1.6ex}[0pt]{$y^2=x^3+x+2$} &
\raisebox{1.6ex}[0pt]{$1\pm2\sqrt{-1}$}&
\raisebox{1.6ex}[0pt]{$\frac{\mathbb{Z}}{4\mathbb{Z}}$}&
\raisebox{1.6ex}[0pt]{$\frac{\mathbb{Z}}{116\mathbb{Z}}$}&
\multicolumn{1}{|c|}{$\lambda(29)=1$}&\multicolumn{1}{|c|}{$K_2(29^m)(29)\cong
\frac{\mathbb{Z}}{29^{m+1}\mathbb{Z}},\;m\geq 0$}\\
\hline

{}&{}&{}&{}&{}&
\multicolumn{1}{|c|}{$\lambda(2)=2$}&\multicolumn{1}{|c|}{$K_2(2^m)(2)\cong
\frac{\mathbb{Z}}{2^{m+1}\mathbb{Z}}\times\frac{\mathbb{Z}}{2^{m+2}\mathbb{Z}},\;m\geq 1$}\\
\cline{6-7} \raisebox{1.6ex}[0pt]{$4$}&
\raisebox{1.6ex}[0pt]{$y^2=x^3+x$} &
\raisebox{1.6ex}[0pt]{$1\pm2\sqrt{-1}$}&
\raisebox{1.6ex}[0pt]{$\frac{\mathbb{Z}}{2\mathbb{Z}}\times\frac{\mathbb{Z}}{2\mathbb{Z}}$}&
\raisebox{1.6ex}[0pt]{$\frac{\mathbb{Z}}{2\mathbb{Z}}\times\frac{\mathbb{Z}}{58\mathbb{Z}}$}&
\multicolumn{1}{|c|}{$\lambda(29)=1$}&\multicolumn{1}{|c|}{$K_2(29^m)(29)\cong
\frac{\mathbb{Z}}{29^{m+1}\mathbb{Z}},\;m\geq 0$}\\
\hline

$5$& {$ y^2=x^3-2x+2$} & $\frac{1\pm\sqrt{-19}}{2}$&
$\frac{\mathbb{Z}}{5\mathbb{Z}}$&$\frac{\mathbb{Z}}{121\mathbb{Z}}$&
$\lambda(11)=1$&{
$K_2(11^m)(11)\cong \frac{\mathbb{Z}}{11^{m+2}\mathbb{Z}},\;m \geq 0$}\\
\hline

{}&{}&{}&{}&{}&
\multicolumn{1}{|c|}{$\lambda(2)=2$}&\multicolumn{1}{|c|}{$K_2(2^m)(2)\cong
\frac{\mathbb{Z}}{2^{m+1}\mathbb{Z}}\times\frac{\mathbb{Z}}{2^{m+1}\mathbb{Z}},\;m\geq 2$}\\
\cline{6-7} {}&{}&{}&{}&{}&
\multicolumn{1}{|c|}{$\lambda(3)=1$}&\multicolumn{1}{|c|}{$K_2(3^m)(3)\cong
\frac{\mathbb{Z}}{3^{m+2}\mathbb{Z}},\;m\geq 0$}\\
\cline{6-7} \raisebox{2.6ex}[0pt]{$6$}&
\raisebox{2.6ex}[0pt]{$y^2=x^3+1$} &
\raisebox{2.6ex}[0pt]{$\pm\sqrt{-5}$}&
\raisebox{2.6ex}[0pt]{$\frac{\mathbb{Z}}{6\mathbb{Z}}$}&
\raisebox{2.6ex}[0pt]{$\frac{\mathbb{Z}}{126\mathbb{Z}}$}&
\multicolumn{1}{|c|}{$\lambda(7)=1$}&\multicolumn{1}{|c|}{$K_2(7^m)(7)\cong
\frac{\mathbb{Z}}{7^{m+1}\mathbb{Z}},\;m\geq 0$}\\
\hline

{}&{}&{}&{}&{}&
\multicolumn{1}{|c|}{$\lambda(2)=2$}&\multicolumn{1}{|c|}{$K_2(2^m)(2)\cong
\frac{\mathbb{Z}}{2^{m+1}\mathbb{Z}}\times\frac{\mathbb{Z}}{2^{m+1}\mathbb{Z}},\;m\geq 2$}\\
\cline{6-7} {}&{}&{}&{}&{}&
\multicolumn{1}{|c|}{$\lambda(3)=1$}&\multicolumn{1}{|c|}{$K_2(3^m)(3)\cong
\frac{\mathbb{Z}}{3^{m+2}\mathbb{Z}},\;m\geq 1$}\\
\cline{6-7} \raisebox{2.6ex}[0pt]{$7$}&
\raisebox{2.6ex}[0pt]{$y^2=x^3+2$} &
\raisebox{2.6ex}[0pt]{$\pm\sqrt{-5}$}&
\raisebox{2.6ex}[0pt]{$\frac{\mathbb{Z}}{6\mathbb{Z}}$}&
\raisebox{2.6ex}[0pt]{$\frac{\mathbb{Z}}{126\mathbb{Z}}$}&
\multicolumn{1}{|c|}{$\lambda(7)=1$}&\multicolumn{1}{|c|}{$K_2(7^m)(7)\cong
\frac{\mathbb{Z}}{7^{m+1}\mathbb{Z}},\;m\geq 0$}\\
\hline

$8$& {$ y^2=x^3+2x+1$} & $\frac{-1\pm\sqrt{-19}}{2}$&
$\frac{\mathbb{Z}}{7\mathbb{Z}}$&$\frac{\mathbb{Z}}{131\mathbb{Z}}$&
$\lambda(131)=1$&{
$K_2(131^m)(131)\cong \frac{\mathbb{Z}}{131^{m+1}\mathbb{Z}},\;m \geq 0$}\\
\hline

{}&{}&{}&{}&{}&
\multicolumn{1}{|c|}{$\lambda(2)=2$}&\multicolumn{1}{|c|}{$K_2(2^m)(2)\cong
\frac{\mathbb{Z}}{2^{m+1}\mathbb{Z}}\times\frac{\mathbb{Z}}{2^{m+2}\mathbb{Z}},\;m\geq 1$}\\
\cline{6-7}
 \raisebox{1.6ex}[0pt]{$9$}&
\raisebox{1.6ex}[0pt]{$y^2=x^3-x$} &
\raisebox{1.6ex}[0pt]{$-1\pm2\sqrt{-1}$}&
\raisebox{1.6ex}[0pt]{$\frac{\mathbb{Z}}{2\mathbb{Z}}\times\frac{\mathbb{Z}}{4\mathbb{Z}}$}&
\raisebox{1.6ex}[0pt]{$\frac{\mathbb{Z}}{4\mathbb{Z}}\times\frac{\mathbb{Z}}{34\mathbb{Z}}$}&
\multicolumn{1}{|c|}{$\lambda(17)=1$}&\multicolumn{1}{|c|}{$K_2(17^m)(17)\cong
\frac{\mathbb{Z}}{17^{m+1}\mathbb{Z}},\;m\geq 0$}\\
\hline

{}&{}&{}&{}&{}&
\multicolumn{1}{|c|}{$\lambda(2)=2$}&\multicolumn{1}{|c|}{$K_2(2^m)(2)\cong
\frac{\mathbb{Z}}{2^{m}\mathbb{Z}}\times\frac{\mathbb{Z}}{2^{m+3}\mathbb{Z}},\;m\geq 1$}\\
\cline{6-7}
 \raisebox{1.6ex}[0pt]{$10$}&
\raisebox{1.6ex}[0pt]{$y^2=x^3-x-1$} &
\raisebox{1.6ex}[0pt]{$-1\pm2\sqrt{-1}$}&
\raisebox{1.6ex}[0pt]{$\frac{\mathbb{Z}}{8\mathbb{Z}}$}&
\raisebox{1.6ex}[0pt]{$\frac{\mathbb{Z}}{136\mathbb{Z}}$}&
\multicolumn{1}{|c|}{$\lambda(17)=1$}&\multicolumn{1}{|c|}{$K_2(17^m)(17)\cong
\frac{\mathbb{Z}}{17^{m+1}\mathbb{Z}},\;m\geq 0$}\\
\hline

{}&{}&{}&{}&{}&
\multicolumn{1}{|c|}{$\lambda(3)=1$}&\multicolumn{1}{|c|}{$K_2(3^m)(3)\cong
\frac{\mathbb{Z}}{3^{m+1}\mathbb{Z}},\;m\geq 0$}\\
\cline{6-7} \raisebox{1.6ex}[0pt]{$11$}&
\raisebox{1.6ex}[0pt]{$y^2=x^3+x+1$} &
\raisebox{1.6ex}[0pt]{$\frac{-3\pm\sqrt{-11}}{2}$}&
\raisebox{1.6ex}[0pt]{$\frac{\mathbb{Z}}{9\mathbb{Z}}$}&
\raisebox{1.6ex}[0pt]{$\frac{\mathbb{Z}}{141\mathbb{Z}}$}&
\multicolumn{1}{|c|}{$\lambda(47)=1$}&\multicolumn{1}{|c|}{$K_2(47^m)(47)\cong
\frac{\mathbb{Z}}{47^{m+1}\mathbb{Z}},\;m\geq 0$}\\
\hline

{}&{}&{}&{}&{}&
\multicolumn{1}{|c|}{$\lambda(2)=2$}&\multicolumn{1}{|c|}{$K_2(2^m)(2)\cong
\frac{\mathbb{Z}}{2^{m+2}\mathbb{Z}}\times\frac{\mathbb{Z}}{2^{m+3}\mathbb{Z}},\;m\geq 1$}\\
 \cline{6-7}\raisebox{1.6ex}[0pt]{$12$}&
\raisebox{1.6ex}[0pt]{$y^2=x^3-2x$} &
\raisebox{1.6ex}[0pt]{$-2\pm\sqrt{-1}$}&
\raisebox{1.6ex}[0pt]{$\frac{\mathbb{Z}}{10\mathbb{Z}}$}&
\raisebox{1.6ex}[0pt]{$\frac{\mathbb{Z}}{146\mathbb{Z}}$}&
\multicolumn{1}{|c|}{$\lambda(73)=1$}&\multicolumn{1}{|c|}{$K_2(73^m)(73)\cong
\frac{\mathbb{Z}}{73^{m+1}\mathbb{Z}},\;m\geq 0$}\\
\hline

\multicolumn{7}{|c|}{\text{The\; elliptic\; curves\; of\; numbers}\;
$n$\;\text{ and}\; $13-n$\;\text{ are\; twisted.}}\\\hline
\end{tabular}
\end{lrbox}
\resizebox{5.5 in}{3.5 in}{\usebox{\tablebox}}

\begin{lrbox}{\tablebox}
\begin{tabular}{|c|c|c|c|c|c|c|c|}
\multicolumn{8}{c}{\bf TABLE\;\;\; II(2)}\\
\multicolumn{8}{c}{(\text{of\; all\; elliptic\; curves}\; E
\;\text{defined\; over}\;
$\mathbb{F}=\mathbb{Z}/5\mathbb{Z}$)}\\

\multicolumn{8}{c}{}\\
\hline

No. &\multicolumn{1}{|c|}{$E$} &\multicolumn{1}{|c|}{$K_2(E)$}&
\multicolumn{1}{|c|}{$K_4(E)$}& \multicolumn{1}{|c|}{$K_6(E)$}&
\multicolumn{1}{|c|}{$K_8(E)$}& \multicolumn{1}{|c|}{$K_{10}(E)$}&
\multicolumn{1}{|c|}{$K_{12}(E)$}\\
\hline

$1$ & $y^2=x^3+2x$& $\frac{\mathbb{Z}}{106\mathbb{Z}}$&
$\frac{\mathbb{Z}}{3026\mathbb{Z}}$&
$\frac{\mathbb{Z}}{77626\mathbb{Z}}$&
$\frac{\mathbb{Z}}{1950626\mathbb{Z}}$&
$\frac{\mathbb{Z}}{48815626\mathbb{Z}}$&
$\frac{\mathbb{Z}}{1220640626\mathbb{Z}}$\\
\hline

$2$ & $y^2=x^3-x+2$& $\frac{\mathbb{Z}}{111\mathbb{Z}}$&
$\frac{\mathbb{Z}}{3051\mathbb{Z}}$&
$\frac{\mathbb{Z}}{77751\mathbb{Z}}$&
$\frac{\mathbb{Z}}{1951251\mathbb{Z}}$&
$\frac{\mathbb{Z}}{48818751\mathbb{Z}}$&
$\frac{\mathbb{Z}}{1220656251\mathbb{Z}}$\\
\hline

$3$ & $y^2=x^3+x+2$& $\frac{\mathbb{Z}}{116\mathbb{Z}}$&
$\frac{\mathbb{Z}}{3076\mathbb{Z}}$&
$\frac{\mathbb{Z}}{77876\mathbb{Z}}$&
$\frac{\mathbb{Z}}{1951876\mathbb{Z}}$&
$\frac{\mathbb{Z}}{48821876\mathbb{Z}}$&
$\frac{\mathbb{Z}}{1220671876\mathbb{Z}}$\\
\hline

$4$ & $y^2=x^3+x$&
$\frac{\mathbb{Z}}{2\mathbb{Z}}\times\frac{\mathbb{Z}}{58\mathbb{Z}}$&
$\frac{\mathbb{Z}}{2\mathbb{Z}}\times\frac{\mathbb{Z}}{1538\mathbb{Z}}$&
$\frac{\mathbb{Z}}{2\mathbb{Z}}\times\frac{\mathbb{Z}}{38938\mathbb{Z}}$&
$\frac{\mathbb{Z}}{2\mathbb{Z}}\times\frac{\mathbb{Z}}{975938\mathbb{Z}}$&
$\frac{\mathbb{Z}}{2\mathbb{Z}}\times\frac{\mathbb{Z}}{24410938\mathbb{Z}}$&
$\frac{\mathbb{Z}}{2\mathbb{Z}}\times\frac{\mathbb{Z}}{610335938\mathbb{Z}}$\\
\hline

$5$ & $y^2=x^3-2x+2$& $\frac{\mathbb{Z}}{121\mathbb{Z}}$&
$\frac{\mathbb{Z}}{3101\mathbb{Z}}$&
$\frac{\mathbb{Z}}{78001\mathbb{Z}}$&
$\frac{\mathbb{Z}}{1952501\mathbb{Z}}$&
$\frac{\mathbb{Z}}{48825001\mathbb{Z}}$&
$\frac{\mathbb{Z}}{1220687501\mathbb{Z}}$\\
\hline

$6$ & $y^2=x^3+1$& $\frac{\mathbb{Z}}{126\mathbb{Z}}$&
$\frac{\mathbb{Z}}{3126\mathbb{Z}}$&
$\frac{\mathbb{Z}}{78126\mathbb{Z}}$&
$\frac{\mathbb{Z}}{1953126\mathbb{Z}}$&
$\frac{\mathbb{Z}}{48828126\mathbb{Z}}$&
$\frac{\mathbb{Z}}{1220703126\mathbb{Z}}$\\
\hline

$7$ & $y^2=x^3+2$& $\frac{\mathbb{Z}}{126\mathbb{Z}}$&
$\frac{\mathbb{Z}}{3126\mathbb{Z}}$&
$\frac{\mathbb{Z}}{78126\mathbb{Z}}$&
$\frac{\mathbb{Z}}{1953126\mathbb{Z}}$&
$\frac{\mathbb{Z}}{48828126\mathbb{Z}}$&
$\frac{\mathbb{Z}}{1220703126\mathbb{Z}}$\\
\hline

$8$ & $y^2=x^3+2x+1$& $\frac{\mathbb{Z}}{131\mathbb{Z}}$&
$\frac{\mathbb{Z}}{3151\mathbb{Z}}$&
$\frac{\mathbb{Z}}{78251\mathbb{Z}}$&
$\frac{\mathbb{Z}}{19537\mathbb{Z}}$&
$\frac{\mathbb{Z}}{48831251\mathbb{Z}}$&
$\frac{\mathbb{Z}}{1220718751\mathbb{Z}}$\\
\hline

$9$ & $y^2=x^3-x$&
$\frac{\mathbb{Z}}{4\mathbb{Z}}\times\frac{\mathbb{Z}}{34\mathbb{Z}}$&
$\frac{\mathbb{Z}}{4\mathbb{Z}}\times\frac{\mathbb{Z}}{794\mathbb{Z}}$&
$\frac{\mathbb{Z}}{4\mathbb{Z}}\times\frac{\mathbb{Z}}{19594\mathbb{Z}}$&
$\frac{\mathbb{Z}}{4\mathbb{Z}}\times\frac{\mathbb{Z}}{488594\mathbb{Z}}$&
$\frac{\mathbb{Z}}{4\mathbb{Z}}\times\frac{\mathbb{Z}}{12208594\mathbb{Z}}$&
$\frac{\mathbb{Z}}{4\mathbb{Z}}\times\frac{\mathbb{Z}}{305183594\mathbb{Z}}$\\
\hline

$10$ & $y^2=x^3-x-1$& $\frac{\mathbb{Z}}{136\mathbb{Z}}$&
$\frac{\mathbb{Z}}{3176\mathbb{Z}}$&
$\frac{\mathbb{Z}}{78376\mathbb{Z}}$&
$\frac{\mathbb{Z}}{1954376\mathbb{Z}}$&
$\frac{\mathbb{Z}}{48834376\mathbb{Z}}$&
$\frac{\mathbb{Z}}{1220734376\mathbb{Z}}$\\
\hline

$11$ & $y^2=x^3+x+1$& $\frac{\mathbb{Z}}{141\mathbb{Z}}$&
$\frac{\mathbb{Z}}{3201\mathbb{Z}}$&
$\frac{\mathbb{Z}}{78501\mathbb{Z}}$&
$\frac{\mathbb{Z}}{1955001\mathbb{Z}}$&
$\frac{\mathbb{Z}}{48837501\mathbb{Z}}$&
$\frac{\mathbb{Z}}{1220750001\mathbb{Z}}$\\
\hline

$12$ & $y^2=x^3-2x$& $\frac{\mathbb{Z}}{146\mathbb{Z}}$&
$\frac{\mathbb{Z}}{3226\mathbb{Z}}$&
$\frac{\mathbb{Z}}{78626\mathbb{Z}}$&
$\frac{\mathbb{Z}}{1955626\mathbb{Z}}$&
$\frac{\mathbb{Z}}{48840626\mathbb{Z}}$&
$\frac{\mathbb{Z}}{1220765626\mathbb{Z}}$\\

\hline \multicolumn{8}{|c|}{\text{The\; elliptic\; curves\; of\;
numbers}\; $n$\;\text{ and}\; $13-n$\; \text{ are\; twisted}.}\\
\hline

\end{tabular}
\end{lrbox}
\resizebox{5.5 in}{3.5 in}{\usebox{\tablebox}}

\begin{lrbox}{\tablebox}
\begin{tabular}{|c|c|c|c|c|c|l|}
\multicolumn{7}{c}{\bf TABLE\;\;\; III(1)}\\
\multicolumn{7}{c}{(\text{of\; all\; elliptic\; curves}\; E
\;\text{defined\; over}\;
$\mathbb{F}=\mathbb{Z}/7\mathbb{Z}$)}\\

\multicolumn{7}{c}{}\\
\hline

 No. &\multicolumn{1}{|c|}{$E$} & \multicolumn{1}{|c|}{\text{roots}}&
\multicolumn{1}{|c|}{$E(\mathbb{F})$}&
\multicolumn{1}{|c|}{$K_2(E)$}&
\multicolumn{1}{|c|}{$\lambda(l)$}& \multicolumn{1}{|c|}{\small $l$-{\text{Sylow}}\;\;\text{subgroups}}\\
\hline

{}&{}&{}&{}&{}&
\multicolumn{1}{|c|}{$\lambda(3)=2$}&\multicolumn{1}{|c|}{$K_2(3^m)(3)\cong
\frac{\mathbb{Z}}{3^{m+1}\mathbb{Z}}\times\frac{\mathbb{Z}}{3^{m+2}\mathbb{Z}},\;m\geq 1$}\\
 \cline{6-7}\raisebox{1.6ex}[0pt]{$1$}&
\raisebox{1.6ex}[0pt]{$y^2=x^3+4$} &
\raisebox{1.6ex}[0pt]{$\frac{5\pm\sqrt{-3}}{2}$}&
\raisebox{1.6ex}[0pt]{$\frac{\mathbb{Z}}{3\mathbb{Z}}$}&
\raisebox{1.6ex}[0pt]{$\frac{\mathbb{Z}}{309\mathbb{Z}}$}&
\multicolumn{1}{|c|}{$\lambda(103)=1$}&\multicolumn{1}{|c|}{$K_2(103^m)(103)\cong
\frac{\mathbb{Z}}{103^{m+1}\mathbb{Z}},\;m\geq 0$}\\

\hline {}&{}&{}&{}&{}&
\multicolumn{1}{|c|}{$\lambda(2)=2$}&\multicolumn{1}{|c|}{$K_2(2^m)(2)\cong
\frac{\mathbb{Z}}{2^{m+1}\mathbb{Z}}\times\frac{\mathbb{Z}}{2^{m+1}\mathbb{Z}},\;m\geq 1$}\\
\cline{6-7} \raisebox{1.6ex}[0pt]{$2$}&
\raisebox{1.6ex}[0pt]{$y^2=x^3+6$} &
\raisebox{1.6ex}[0pt]{$2\pm\sqrt{-3}$}&
\raisebox{1.6ex}[0pt]{$\frac{\mathbb{Z}}{2\mathbb{Z}}\times\frac{\mathbb{Z}}{2\mathbb{Z}}$}&
\raisebox{1.6ex}[0pt]{$\frac{\mathbb{Z}}{2\mathbb{Z}}\times\frac{\mathbb{Z}}{158\mathbb{Z}}$}&
\multicolumn{1}{|c|}{$\lambda(79)=1$}&\multicolumn{1}{|c|}{$K_2(79^m)(79)\cong
\frac{\mathbb{Z}}{79^{m+1}\mathbb{Z}},\;m\geq 0$}\\

\hline {}&{}&{}&{}&{}&
\multicolumn{1}{|c|}{$\lambda(2)=2$}&\multicolumn{1}{|c|}{$K_2(2^m)(2)\cong
\frac{\mathbb{Z}}{2^{m+1}\mathbb{Z}}\times\frac{\mathbb{Z}}{2^{m+1}\mathbb{Z}},\;m\geq 1$}\\
 \cline{6-7}\raisebox{1.6ex}[0pt]{$3$}&
\raisebox{1.6ex}[0pt]{$y^2=x^3+6x+6$} &
\raisebox{1.6ex}[0pt]{$2\pm\sqrt{-3}$}&
\raisebox{1.6ex}[0pt]{$\frac{\mathbb{Z}}{4\mathbb{Z}}$}&
\raisebox{1.6ex}[0pt]{$\frac{\mathbb{Z}}{316\mathbb{Z}}$}&
\multicolumn{1}{|c|}{$\lambda(79)=1$}&\multicolumn{1}{|c|}{$K_2(79^m)(79)\cong
\frac{\mathbb{Z}}{79^{m+1}\mathbb{Z}},\;m\geq 0$}\\

\hline {}&{}&{}&{}&{}&
\multicolumn{1}{|c|}{$\lambda(17)=1$}&\multicolumn{1}{|c|}{$K_2(17^m)(17)\cong
\frac{\mathbb{Z}}{17^{m+1}\mathbb{Z}},\;m\geq 0$}\\
\cline{6-7} \raisebox{1.6ex}[0pt]{$4$}&
\raisebox{1.6ex}[0pt]{$y^2=x^3+x+1$} &
\raisebox{1.6ex}[0pt]{$\frac{3\pm\sqrt{-19}}{2}$}&
\raisebox{1.6ex}[0pt]{$\frac{\mathbb{Z}}{5\mathbb{Z}}$}&
\raisebox{1.6ex}[0pt]{$\frac{\mathbb{Z}}{323\mathbb{Z}}$}&
\multicolumn{1}{|c|}{$\lambda(19)=2$}&\multicolumn{1}{|c|}{$K_2(19^m)(19)\cong
\frac{\mathbb{Z}}{19^{m}\mathbb{Z}}\times\frac{\mathbb{Z}}{19^{m+1}\mathbb{Z}},\;m\geq 1$}\\

\hline {}&{}&{}&{}&{}&
\multicolumn{1}{|c|}{$\lambda(2)=2$}&\multicolumn{1}{|c|}{$K_2(2^m)(2)\cong
\frac{\mathbb{Z}}{2^{m}\mathbb{Z}}\times\frac{\mathbb{Z}}{2^{m+1}\mathbb{Z}},\;m\geq 2$}\\
\cline{6-7}
 {}&{}&{}&{}&{}&
\multicolumn{1}{|c|}{$\lambda(3)=2$}&\multicolumn{1}{|c|}{$K_2(3^m)(3)\cong
\frac{\mathbb{Z}}{3^{m}\mathbb{Z}}\times\frac{\mathbb{Z}}{3^{m+1}\mathbb{Z}},\;m\geq 1$}\\
\cline{6-7} {}&{}&{}&{}&{}&
\multicolumn{1}{|c|}{$\lambda(5)=1$}&\multicolumn{1}{|c|}{$K_2(5^m)(5)\cong
\frac{\mathbb{Z}}{5^{m+1}\mathbb{Z}},\;m\geq 0$}\\
\cline{6-7} \raisebox{4.5ex}[0pt]{$5$}&
\raisebox{4.5ex}[0pt]{$y^2=x^3+3x+3$} &
\raisebox{4.5ex}[0pt]{$1\pm\sqrt{-6}$}&
\raisebox{4.5ex}[0pt]{$\frac{\mathbb{Z}}{6\mathbb{Z}}$}&
\raisebox{4.5ex}[0pt]{$\frac{\mathbb{Z}}{330\mathbb{Z}}$}&
\multicolumn{1}{|c|}{$\lambda(11)=1$}&\multicolumn{1}{|c|}{$K_2(11^m)(11)\cong
\frac{\mathbb{Z}}{11^{m+1}\mathbb{Z}},\;m\geq 0$}\\

\hline {}&{}&{}&{}&{}&
\multicolumn{1}{|c|}{$\lambda(2)=2$}&\multicolumn{1}{|c|}{$K_2(2^m)(2)\cong
\frac{\mathbb{Z}}{2^{m}\mathbb{Z}}\times\frac{\mathbb{Z}}{2^{m+1}\mathbb{Z}},\;m\geq 2$}\\
\cline{6-7}
 {}&{}&{}&{}&{}&
\multicolumn{1}{|c|}{$\lambda(3)=2$}&\multicolumn{1}{|c|}{$K_2(3^m)(3)\cong
\frac{\mathbb{Z}}{3^{m}\mathbb{Z}}\times\frac{\mathbb{Z}}{3^{m+1}\mathbb{Z}},\;m\geq 1$}\\
\cline{6-7} {}&{}&{}&{}&{}&
\multicolumn{1}{|c|}{$\lambda(5)=1$}&\multicolumn{1}{|c|}{$K_2(5^m)(5)\cong
\frac{\mathbb{Z}}{5^{m+1}\mathbb{Z}},\;m\geq 0$}\\
\cline{6-7} \raisebox{4.5ex}[0pt]{$6$}&
\raisebox{4.5ex}[0pt]{$y^2=x^3+x+3$} &
\raisebox{4.5ex}[0pt]{$1\pm\sqrt{-6}$}&
\raisebox{4.5ex}[0pt]{$\frac{\mathbb{Z}}{6\mathbb{Z}}$}&
\raisebox{4.5ex}[0pt]{$\frac{\mathbb{Z}}{330\mathbb{Z}}$}&
\multicolumn{1}{|c|}{$\lambda(11)=1$}&\multicolumn{1}{|c|}{$K_2(11^m)(11)\cong
\frac{\mathbb{Z}}{11^{m+1}\mathbb{Z}},\;m\geq 0$}\\

\hline $7$&
$y^2=x^3+5$&$\frac{1\pm3\sqrt{-3}}{2}$&$\frac{\mathbb{Z}}{7\mathbb{Z}}$&
$\frac{\mathbb{Z}}{337\mathbb{Z}}$&$\lambda(337)=1$&$K_2(337^m)(337)\cong
\frac{\mathbb{Z}}{337^{m+1}\mathbb{Z}},\;m\geq 0$\\

\hline $8$&
$y^2=x^3+6x+5$&$\frac{1\pm3\sqrt{-3}}{2}$&$\frac{\mathbb{Z}}{7\mathbb{Z}}$&
$\frac{\mathbb{Z}}{337\mathbb{Z}}$&$\lambda(337)=1$&$K_2(337^m)(337)\cong
\frac{\mathbb{Z}}{337^{m+1}\mathbb{Z}},\;m\geq 0$\\

\hline {}&{}&{}&{}&{}&
\multicolumn{1}{|c|}{$\lambda(2)=2$}&\multicolumn{1}{|c|}{$K_2(2^m)(2)\cong
\frac{\mathbb{Z}}{2^{m+2}\mathbb{Z}}\times\frac{\mathbb{Z}}{2^{m+2}\mathbb{Z}},\;m\geq 1$}\\
\cline{6-7}
 \raisebox{1.6ex}[0pt]{$9$}&
\raisebox{1.6ex}[0pt]{$y^2=x^3+x$} &
\raisebox{1.6ex}[0pt]{$\pm\sqrt{-7}$}&
\raisebox{1.6ex}[0pt]{$\frac{\mathbb{Z}}{8\mathbb{Z}}$}&
\raisebox{1.6ex}[0pt]{$\frac{\mathbb{Z}}{344\mathbb{Z}}$}&
\multicolumn{1}{|c|}{$\lambda(43)=1$}&\multicolumn{1}{|c|}{$K_2(43^m)(43)\cong
\frac{\mathbb{Z}}{43^{m+1}\mathbb{Z}},\;m\geq 0$}\\

\hline {}&{}&{}&{}&{}&
\multicolumn{1}{|c|}{$\lambda(2)=2$}&\multicolumn{1}{|c|}{$K_2(2^m)(2)\cong
\frac{\mathbb{Z}}{2^{m+2}\mathbb{Z}}\times\frac{\mathbb{Z}}{2^{m+2}\mathbb{Z}},\;m\geq 1$}\\
\cline{6-7}
 \raisebox{1.6ex}[0pt]{$10$}&
\raisebox{1.6ex}[0pt]{$y^2=x^3+3x$} &
\raisebox{1.6ex}[0pt]{$\pm\sqrt{-7}$}&
\raisebox{1.6ex}[0pt]{$\frac{\mathbb{Z}}{2\mathbb{Z}}\times\frac{\mathbb{Z}}{4\mathbb{Z}}$}&
\raisebox{1.6ex}[0pt]{$\frac{\mathbb{Z}}{4\mathbb{Z}}\times\frac{\mathbb{Z}}{86\mathbb{Z}}$}&
\multicolumn{1}{|c|}{$\lambda(43)=1$}&\multicolumn{1}{|c|}{$K_2(43^m)(43)\cong
\frac{\mathbb{Z}}{43^{m+1}\mathbb{Z}},\;m\geq 0$}\\

\hline {}&{}&{}&{}&{}&
\multicolumn{1}{|c|}{$\lambda(3)=2$}&\multicolumn{1}{|c|}{$K_2(3^m)(3)\cong
\frac{\mathbb{Z}}{3^{m}\mathbb{Z}}\times\frac{\mathbb{Z}}{3^{m+3}\mathbb{Z}},\;m\geq 1$}\\
\cline{6-7}
 \raisebox{1.6ex}[0pt]{$11$}&
\raisebox{1.6ex}[0pt]{$y^2=x^3+3x+2$} &
\raisebox{1.6ex}[0pt]{$\frac{-1\pm3\sqrt{-3}}{2}$}&
\raisebox{1.6ex}[0pt]{$\frac{\mathbb{Z}}{9\mathbb{Z}}$}&
\raisebox{1.6ex}[0pt]{$\frac{\mathbb{Z}}{351\mathbb{Z}}$}&
\multicolumn{1}{|c|}{$\lambda(13)=1$}&\multicolumn{1}{|c|}{$K_2(13^m)(13)\cong
\frac{\mathbb{Z}}{13^{m+1}\mathbb{Z}},\;m\geq 0$}\\

\hline {}&{}&{}&{}&{}&
\multicolumn{1}{|c|}{$\lambda(3)=2$}&\multicolumn{1}{|c|}{$K_2(3^m)(3)\cong
\frac{\mathbb{Z}}{3^{m+1}\mathbb{Z}}\times\frac{\mathbb{Z}}{3^{m+2}\mathbb{Z}},\;m\geq 1$}\\
\cline{6-7}
 \raisebox{1.6ex}[0pt]{$12$}&
\raisebox{1.6ex}[0pt]{$y^2=x^3+2$} &
\raisebox{1.6ex}[0pt]{$\frac{-1\pm3\sqrt{-3}}{2}$}&
\raisebox{1.6ex}[0pt]{$\frac{\mathbb{Z}}{3\mathbb{Z}}\times\frac{\mathbb{Z}}{3\mathbb{Z}}$}&
\raisebox{1.6ex}[0pt]{$\frac{\mathbb{Z}}{9\mathbb{Z}}\times\frac{\mathbb{Z}}{39\mathbb{Z}}$}&
\multicolumn{1}{|c|}{$\lambda(13)=1$}&\multicolumn{1}{|c|}{$K_2(13^m)(13)\cong
\frac{\mathbb{Z}}{13^{m+1}\mathbb{Z}},\;m\geq 0$}\\

\hline {}&{}&{}&{}&{}&
\multicolumn{1}{|c|}{$\lambda(2)=2$}&\multicolumn{1}{|c|}{$K_2(2^m)(2)\cong
\frac{\mathbb{Z}}{2^{m}\mathbb{Z}}\times\frac{\mathbb{Z}}{2^{m+1}\mathbb{Z}},\;m\geq 2$}\\
\cline{6-7}
 \raisebox{1.6ex}[0pt]{$13$}&
\raisebox{1.6ex}[0pt]{$y^2=x^3+2x+4$} &
\raisebox{1.6ex}[0pt]{$-1\pm\sqrt{-6}$}&
\raisebox{1.6ex}[0pt]{$\frac{\mathbb{Z}}{10\mathbb{Z}}$}&
\raisebox{1.6ex}[0pt]{$\frac{\mathbb{Z}}{358\mathbb{Z}}$}&
\multicolumn{1}{|c|}{$\lambda(179)=1$}&\multicolumn{1}{|c|}{$K_2(179^m)(179)\cong
\frac{\mathbb{Z}}{179^{m+1}\mathbb{Z}},\;m\geq 0$}\\

\hline {}&{}&{}&{}&{}&
\multicolumn{1}{|c|}{$\lambda(2)=2$}&\multicolumn{1}{|c|}{$K_2(2^m)(2)\cong
\frac{\mathbb{Z}}{2^{m}\mathbb{Z}}\times\frac{\mathbb{Z}}{2^{m+1}\mathbb{Z}},\;m\geq 2$}\\
\cline{6-7}
 \raisebox{1.6ex}[0pt]{$14$}&
\raisebox{1.6ex}[0pt]{$y^2=x^3+6x+4$} &
\raisebox{1.6ex}[0pt]{$-1\pm\sqrt{-6}$}&
\raisebox{1.6ex}[0pt]{$\frac{\mathbb{Z}}{10\mathbb{Z}}$}&
\raisebox{1.6ex}[0pt]{$\frac{\mathbb{Z}}{358\mathbb{Z}}$}&
\multicolumn{1}{|c|}{$\lambda(179)=1$}&\multicolumn{1}{|c|}{$K_2(179^m)(179)\cong
\frac{\mathbb{Z}}{179^{m+1}\mathbb{Z}},\;m\geq 0$}\\

\hline {}&{}&{}&{}&{}&
\multicolumn{1}{|c|}{$\lambda(5)=1$}&\multicolumn{1}{|c|}{$K_2(5^m)(5)\cong
\frac{\mathbb{Z}}{5^{m+1}\mathbb{Z}},\;m\geq 0$}\\
\cline{6-7}
 \raisebox{1.6ex}[0pt]{$15$}&
\raisebox{1.6ex}[0pt]{$y^2=x^3+2x+6$} &
\raisebox{1.6ex}[0pt]{$\frac{-3\pm\sqrt{-19}}{2}$}&
\raisebox{1.6ex}[0pt]{$\frac{\mathbb{Z}}{11\mathbb{Z}}$}&
\raisebox{1.6ex}[0pt]{$\frac{\mathbb{Z}}{365\mathbb{Z}}$}&
\multicolumn{1}{|c|}{$\lambda(73)=1$}&\multicolumn{1}{|c|}{$K_2(73^m)(73)\cong
\frac{\mathbb{Z}}{73^{m+1}\mathbb{Z}},\;m\geq 0$}\\

\hline {}&{}&{}&{}&{}&
\multicolumn{1}{|c|}{$\lambda(2)=2$}&\multicolumn{1}{|c|}{$K_2(2^m)(2)\cong
\frac{\mathbb{Z}}{2^{m+1}\mathbb{Z}}\times\frac{\mathbb{Z}}{2^{m+1}\mathbb{Z}},\;m\geq 1$}\\
\cline{6-7} {}&{}&{}&{}&{}&
\multicolumn{1}{|c|}{$\lambda(3)=2$}&\multicolumn{1}{|c|}{$K_2(3^m)(3)\cong
\frac{\mathbb{Z}}{3^{m+1}\mathbb{Z}}\times\frac{\mathbb{Z}}{3^{m+2}\mathbb{Z}},\;m\geq 1$}\\
\cline{6-7} \raisebox{2.6ex}[0pt]{$16$}&
\raisebox{2.6ex}[0pt]{$y^2=x^3+3x+1$} &
\raisebox{2.6ex}[0pt]{$-2\pm\sqrt{-3}$}&
\raisebox{2.6ex}[0pt]{$\frac{\mathbb{Z}}{12\mathbb{Z}}$}&
\raisebox{2.6ex}[0pt]{$\frac{\mathbb{Z}}{372\mathbb{Z}}$}&
\multicolumn{1}{|c|}{$\lambda(31)=1$}&\multicolumn{1}{|c|}{$K_2(31^m)(31)\cong
\frac{\mathbb{Z}}{31^{m+1}\mathbb{Z}},\;m\geq 0$}\\

\hline {}&{}&{}&{}&{}&
\multicolumn{1}{|c|}{$\lambda(2)=2$}&\multicolumn{1}{|c|}{$K_2(2^m)(2)\cong
\frac{\mathbb{Z}}{2^{m+1}\mathbb{Z}}\times\frac{\mathbb{Z}}{2^{m+1}\mathbb{Z}},\;m\geq 1$}\\
\cline{6-7} {}&{}&{}&{}&{}&
\multicolumn{1}{|c|}{$\lambda(3)=2$}&\multicolumn{1}{|c|}{$K_2(3^m)(3)\cong
\frac{\mathbb{Z}}{3^{m+1}\mathbb{Z}}\times\frac{\mathbb{Z}}{3^{m+2}\mathbb{Z}},\;m\geq 1$}\\
\cline{6-7} \raisebox{2.6ex}[0pt]{$17$}&
\raisebox{2.6ex}[0pt]{$y^2=x^3+1$} &
\raisebox{2.6ex}[0pt]{$-2\pm\sqrt{-3}$}&
\raisebox{2.6ex}[0pt]{$\frac{\mathbb{Z}}{2\mathbb{Z}}\times\frac{\mathbb{Z}}{6\mathbb{Z}}$}&
\raisebox{2.6ex}[0pt]{$\frac{\mathbb{Z}}{2\mathbb{Z}}\times\frac{\mathbb{Z}}{186\mathbb{Z}}$}&
\multicolumn{1}{|c|}{$\lambda(31)=1$}&\multicolumn{1}{|c|}{$K_2(31^m)(31)\cong
\frac{\mathbb{Z}}{31^{m+1}\mathbb{Z}},\;m\geq 0$}\\

\hline $18$&
$y^2=x^3+3$&$\frac{-5\pm\sqrt{-3}}{2}$&$\frac{\mathbb{Z}}{13\mathbb{Z}}$&
$\frac{\mathbb{Z}}{379\mathbb{Z}}$&$\lambda(379)=1$&$K_2(379^m)(379)\cong
\frac{\mathbb{Z}}{379^{m+1}\mathbb{Z}},\;m\geq 0$\\

\hline \multicolumn{7}{|c|}{\text{The\; elliptic\; curves\; of\;
numbers}\; $n$\;\text{ and}\; $19-n$\; \text{ are\;
twisted}.}\\\hline

\end{tabular}
\end{lrbox}
\resizebox{5.5 in}{!}{\usebox{\tablebox}}

\begin{lrbox}{\tablebox}
\begin{tabular}{|c|c|c|c|c|c|c|c|}
\multicolumn{8}{c}{\bf TABLE\;\;\; III(2)}\\
\multicolumn{8}{c}{(\text{of\; all\; elliptic\; curves}\; E
\;\text{defined\; over}\;
$\mathbb{F}=\mathbb{Z}/7\mathbb{Z}$)}\\

\multicolumn{8}{c}{}\\
\hline

No. &\multicolumn{1}{|c|}{$E$} &\multicolumn{1}{|c|}{$K_2(E)$}&
\multicolumn{1}{|c|}{$K_4(E)$}& \multicolumn{1}{|c|}{$K_6(E)$}&
\multicolumn{1}{|c|}{$K_8(E)$}& \multicolumn{1}{|c|}{$K_{10}(E)$}&
\multicolumn{1}{|c|}{$K_{12}(E)$}\\
\hline

$1$ & $y^2=x^3+4$& $\frac{\mathbb{Z}}{309\mathbb{Z}}$&
$\frac{\mathbb{Z}}{16563\mathbb{Z}}$&
$\frac{\mathbb{Z}}{821829\mathbb{Z}}$&
$\frac{\mathbb{Z}}{40341603\mathbb{Z}}$&
$\frac{\mathbb{Z}}{1977242709\mathbb{Z}}$&
$\frac{\mathbb{Z}}{96888422163\mathbb{Z}}$\\
\hline

$2$ & $y^2=x^3+6$&
$\frac{\mathbb{Z}}{2\mathbb{Z}}\times\frac{\mathbb{Z}}{158\mathbb{Z}}$&
$\frac{\mathbb{Z}}{2\mathbb{Z}}\times\frac{\mathbb{Z}}{8306\mathbb{Z}}$&
$\frac{\mathbb{Z}}{2\mathbb{Z}}\times\frac{\mathbb{Z}}{411086\mathbb{Z}}$&
$\frac{\mathbb{Z}}{2\mathbb{Z}}\times\frac{\mathbb{Z}}{20172002\mathbb{Z}}$&
$\frac{\mathbb{Z}}{2\mathbb{Z}}\times\frac{\mathbb{Z}}{988629758\mathbb{Z}}$&
$\frac{\mathbb{Z}}{2\mathbb{Z}}\times\frac{\mathbb{Z}}{48444269906\mathbb{Z}}$\\
\hline

$3$ & $y^2=x^3+6x+6$& $\frac{\mathbb{Z}}{316\mathbb{Z}}$&
$\frac{\mathbb{Z}}{16612\mathbb{Z}}$&
$\frac{\mathbb{Z}}{822172\mathbb{Z}}$&
$\frac{\mathbb{Z}}{40344004\mathbb{Z}}$&
$\frac{\mathbb{Z}}{1977259516\mathbb{Z}}$&
$\frac{\mathbb{Z}}{96888539812\mathbb{Z}}$\\
\hline

$4$ & $y^2=x^3+x+1$& $\frac{\mathbb{Z}}{323\mathbb{Z}}$&
$\frac{\mathbb{Z}}{16661\mathbb{Z}}$&
$\frac{\mathbb{Z}}{822515\mathbb{Z}}$&
$\frac{\mathbb{Z}}{40346405\mathbb{Z}}$&
$\frac{\mathbb{Z}}{1977276323\mathbb{Z}}$&
$\frac{\mathbb{Z}}{96888657461\mathbb{Z}}$\\
\hline

$5$ & $y^2=x^3+3x+3$& $\frac{\mathbb{Z}}{330\mathbb{Z}}$&
$\frac{\mathbb{Z}}{16710\mathbb{Z}}$&
$\frac{\mathbb{Z}}{822858\mathbb{Z}}$&
$\frac{\mathbb{Z}}{40348806\mathbb{Z}}$&
$\frac{\mathbb{Z}}{1977293130\mathbb{Z}}$&
$\frac{\mathbb{Z}}{96888775110\mathbb{Z}}$\\
\hline

$6$ & $y^2=x^3+x+3$& $\frac{\mathbb{Z}}{330\mathbb{Z}}$&
$\frac{\mathbb{Z}}{16710\mathbb{Z}}$&
$\frac{\mathbb{Z}}{822858\mathbb{Z}}$&
$\frac{\mathbb{Z}}{40348806\mathbb{Z}}$&
$\frac{\mathbb{Z}}{1977293130\mathbb{Z}}$&
$\frac{\mathbb{Z}}{96888775110\mathbb{Z}}$\\
\hline

$7$ & $y^2=x^3+5$& $\frac{\mathbb{Z}}{337\mathbb{Z}}$&
$\frac{\mathbb{Z}}{16759\mathbb{Z}}$&
$\frac{\mathbb{Z}}{823201\mathbb{Z}}$&
$\frac{\mathbb{Z}}{40351207\mathbb{Z}}$&
$\frac{\mathbb{Z}}{1977309937\mathbb{Z}}$&
$\frac{\mathbb{Z}}{96888892759\mathbb{Z}}$\\
\hline

$8$ & $y^2=x^3+6x+5$& $\frac{\mathbb{Z}}{337\mathbb{Z}}$&
$\frac{\mathbb{Z}}{16759\mathbb{Z}}$&
$\frac{\mathbb{Z}}{823201\mathbb{Z}}$&
$\frac{\mathbb{Z}}{40351207\mathbb{Z}}$&
$\frac{\mathbb{Z}}{1977309937\mathbb{Z}}$&
$\frac{\mathbb{Z}}{96888892759\mathbb{Z}}$\\
\hline

$9$ & $y^2=x^3+x$& $\frac{\mathbb{Z}}{334\mathbb{Z}}$&
$\frac{\mathbb{Z}}{16808\mathbb{Z}}$&
$\frac{\mathbb{Z}}{823544\mathbb{Z}}$&
$\frac{\mathbb{Z}}{40353608\mathbb{Z}}$&
$\frac{\mathbb{Z}}{1977326744\mathbb{Z}}$&
$\frac{\mathbb{Z}}{96889010408\mathbb{Z}}$\\
\hline

$10$ & $y^2=x^3+3x$&
$\frac{\mathbb{Z}}{4\mathbb{Z}}\times\frac{\mathbb{Z}}{86\mathbb{Z}}$&
$\frac{\mathbb{Z}}{4\mathbb{Z}}\times\frac{\mathbb{Z}}{4202\mathbb{Z}}$&
$\frac{\mathbb{Z}}{4\mathbb{Z}}\times\frac{\mathbb{Z}}{205886\mathbb{Z}}$&
$\frac{\mathbb{Z}}{4\mathbb{Z}}\times\frac{\mathbb{Z}}{10088402\mathbb{Z}}$&
$\frac{\mathbb{Z}}{4\mathbb{Z}}\times\frac{\mathbb{Z}}{494331686\mathbb{Z}}$&
$\frac{\mathbb{Z}}{4\mathbb{Z}}\times\frac{\mathbb{Z}}{24222252602\mathbb{Z}}$\\
\hline

$11$ & $y^2=x^3+3x+2$& $\frac{\mathbb{Z}}{351\mathbb{Z}}$&
$\frac{\mathbb{Z}}{16857\mathbb{Z}}$&
$\frac{\mathbb{Z}}{823887\mathbb{Z}}$&
$\frac{\mathbb{Z}}{40356009\mathbb{Z}}$&
$\frac{\mathbb{Z}}{1977343551\mathbb{Z}}$&
$\frac{\mathbb{Z}}{96889128057\mathbb{Z}}$\\
\hline

$12$ & $y^2=x^3+2$&
$\frac{\mathbb{Z}}{9\mathbb{Z}}\times\frac{\mathbb{Z}}{39\mathbb{Z}}$&
$\frac{\mathbb{Z}}{3\mathbb{Z}}\times\frac{\mathbb{Z}}{5619\mathbb{Z}}$&
$\frac{\mathbb{Z}}{3\mathbb{Z}}\times\frac{\mathbb{Z}}{274629\mathbb{Z}}$&
$\frac{\mathbb{Z}}{3\mathbb{Z}}\times\frac{\mathbb{Z}}{13452003\mathbb{Z}}$&
$\frac{\mathbb{Z}}{3\mathbb{Z}}\times\frac{\mathbb{Z}}{659114517\mathbb{Z}}$&
$\frac{\mathbb{Z}}{3\mathbb{Z}}\times\frac{\mathbb{Z}}{32296376019\mathbb{Z}}$\\
\hline

$13$ & $y^2=x^3+2x+4$& $\frac{\mathbb{Z}}{358\mathbb{Z}}$&
$\frac{\mathbb{Z}}{16906\mathbb{Z}}$&
$\frac{\mathbb{Z}}{824230\mathbb{Z}}$&
$\frac{\mathbb{Z}}{40358410\mathbb{Z}}$&
$\frac{\mathbb{Z}}{1977360358\mathbb{Z}}$&
$\frac{\mathbb{Z}}{96889245706\mathbb{Z}}$\\
\hline

$14$ & $y^2=x^3+6x+4$& $\frac{\mathbb{Z}}{358\mathbb{Z}}$&
$\frac{\mathbb{Z}}{16906\mathbb{Z}}$&
$\frac{\mathbb{Z}}{824230\mathbb{Z}}$&
$\frac{\mathbb{Z}}{40358410\mathbb{Z}}$&
$\frac{\mathbb{Z}}{1977360358\mathbb{Z}}$&
$\frac{\mathbb{Z}}{96889245706\mathbb{Z}}$\\
\hline

$15$ & $y^2=x^3+2x+6$& $\frac{\mathbb{Z}}{365\mathbb{Z}}$&
$\frac{\mathbb{Z}}{16955\mathbb{Z}}$&
$\frac{\mathbb{Z}}{824573\mathbb{Z}}$&
$\frac{\mathbb{Z}}{40360811\mathbb{Z}}$&
$\frac{\mathbb{Z}}{1977377165\mathbb{Z}}$&
$\frac{\mathbb{Z}}{96889363355\mathbb{Z}}$\\
\hline

$16$ & $y^2=x^3+3x+1$& $\frac{\mathbb{Z}}{372\mathbb{Z}}$&
$\frac{\mathbb{Z}}{17004\mathbb{Z}}$&
$\frac{\mathbb{Z}}{824916\mathbb{Z}}$&
$\frac{\mathbb{Z}}{40363212\mathbb{Z}}$&
$\frac{\mathbb{Z}}{1977393972\mathbb{Z}}$&
$\frac{\mathbb{Z}}{96889481004\mathbb{Z}}$\\
\hline

$17$ & $y^2=x^3+1$&
$\frac{\mathbb{Z}}{2\mathbb{Z}}\times\frac{\mathbb{Z}}{186\mathbb{Z}}$&
$\frac{\mathbb{Z}}{2\mathbb{Z}}\times\frac{\mathbb{Z}}{8502\mathbb{Z}}$&
$\frac{\mathbb{Z}}{2\mathbb{Z}}\times\frac{\mathbb{Z}}{412458\mathbb{Z}}$&
$\frac{\mathbb{Z}}{2\mathbb{Z}}\times\frac{\mathbb{Z}}{20181606\mathbb{Z}}$&
$\frac{\mathbb{Z}}{2\mathbb{Z}}\times\frac{\mathbb{Z}}{988696986\mathbb{Z}}$&
$\frac{\mathbb{Z}}{2\mathbb{Z}}\times\frac{\mathbb{Z}}{48444740502\mathbb{Z}}$\\
\hline

$18$ & $y^2=x^3+3$& $\frac{\mathbb{Z}}{379\mathbb{Z}}$&
$\frac{\mathbb{Z}}{17053\mathbb{Z}}$&
$\frac{\mathbb{Z}}{825259\mathbb{Z}}$&
$\frac{\mathbb{Z}}{40365613\mathbb{Z}}$&
$\frac{\mathbb{Z}}{1977410779\mathbb{Z}}$&
$\frac{\mathbb{Z}}{96889598653\mathbb{Z}}$\\
\hline

\multicolumn{8}{|c|}{\text{The\; elliptic\; curves\; of\; numbers}\;
$n$\;\text{ and}\; $19-n$\; \text{ are\; twisted}.}\\
\hline

\end{tabular}
\end{lrbox}
\resizebox{5.5 in}{3.5 in}{\usebox{\tablebox}}

\vskip 2cm
\begin{lrbox}{\tablebox}
\begin{tabular}{|c|c|c|c|c|c|}
\multicolumn{6}{c}{\bf TABLE\;\;\; IV}\\
\multicolumn{6}{c}{(\text{of\; all\; elliptic\; curves}\; E
\;\text{defined\; over}\;
$\mathbb{F}=\mathbb{Z}/11\mathbb{Z}$)}\\

\multicolumn{6}{c}{}\\
\hline
 No. &\multicolumn{1}{|c|}{$E$} & \multicolumn{1}{|c|}{\text{inverse\;roots}}&
\multicolumn{1}{|c|}{$E(\mathbb{F})$}&
\multicolumn{1}{|c|}{$K_2(E)$}&
\multicolumn{1}{|c|}{$\lambda(l)$}\\

\hline
$1$&$y^2=x^3+4x+2$&$3\pm\sqrt{-2}$&$\frac{\mathbb{Z}}{6\mathbb{Z}}$&
$\frac{\mathbb{Z}}{1266\mathbb{Z}}$&$\lambda(2)=2,\;\lambda(3)=1,\;\lambda(211)=1$\\
\hline

$2$&$y^2=x^3+8x-1$&$\frac{5\pm\sqrt{-19}}{2}$&$\frac{\mathbb{Z}}{7\mathbb{Z}}$&$\frac{\mathbb{Z}}{1277\mathbb{Z}}$&$\lambda(1277)=1$\\
\hline

$3$&$y^2=x^3+4x+5$&$2\pm\sqrt{-7}$&$\frac{\mathbb{Z}}{2\mathbb{Z}}\times\frac{\mathbb{Z}}{4\mathbb{Z}}$&
$\frac{\mathbb{Z}}{2\mathbb{Z}}\times\frac{\mathbb{Z}}{644\mathbb{Z}}$&
$\lambda(2)=\lambda(7)=2,\;\lambda(23)=1$\\

\hline
$4$&$y^2=x^3+8x+8$&$2\pm\sqrt{-7}$&$\frac{\mathbb{Z}}{8\mathbb{Z}}$&
$\frac{\mathbb{Z}}{1288\mathbb{Z}}$&$\lambda(2)=\lambda(7)=2,\;\lambda(23)=1$\\
\hline

$5$&$y^2=x^3+x+4$&$\frac{3\pm\sqrt{-35}}{2}$&$\frac{\mathbb{Z}}{9\mathbb{Z}}$&$\frac{\mathbb{Z}}{1299\mathbb{Z}}$&
$\lambda(3)=\lambda(433)=1$\\
\hline

$6$&$y^2=x^3+2x+2$&$\frac{3\pm\sqrt{-35}}{2}$&$\frac{\mathbb{Z}}{9\mathbb{Z}}$&$\frac{\mathbb{Z}}{1299\mathbb{Z}}$&
$\lambda(3)=\lambda(433)=1$\\
\hline

$7$&$y^2=x^3+4x=8$&$1\pm\sqrt{-10}$&$\frac{\mathbb{Z}}{10\mathbb{Z}}$&
$\frac{\mathbb{Z}}{1310\mathbb{Z}}$&$\lambda(2)=\lambda(5)=2,\;\lambda(131)=1$\\
\hline

$8$&$y^2=x^3+2x+5$&$1\pm\sqrt{-10}$&$\frac{\mathbb{Z}}{10\mathbb{Z}}$&
$\frac{\mathbb{Z}}{1310\mathbb{Z}}$&$\lambda(2)=\lambda(5)=2,\;\lambda(131)=1$\\
\hline

$9$&$y^2=x^3+x+5$&$\frac{1\pm\sqrt{-43}}{2}$&$\frac{\mathbb{Z}}{11\mathbb{Z}}$&
$\frac{\mathbb{Z}}{1321\mathbb{Z}}$&$\lambda(1321)=1$\\
 \hline

$10$&$y^2=x^3+x$&$\pm\sqrt{-11}$&$\frac{\mathbb{Z}}{12\mathbb{Z}}$&
$\frac{\mathbb{Z}}{1332\mathbb{Z}}$&$\lambda(2)=2,\;\lambda(3)=\lambda(37)=1$\\
 \hline

$11$&$y^2=x^3+x$&$\pm\sqrt{-11}$&$\frac{\mathbb{Z}}{12\mathbb{Z}}$&
$\frac{\mathbb{Z}}{1332\mathbb{Z}}$&$\lambda(2)=2,\;\lambda(3)=\lambda(37)=1$\\
 \hline

$12$&$y^2=x^3+2x$&$\pm\sqrt{-11}$&$\frac{\mathbb{Z}}{12\mathbb{Z}}$&
$\frac{\mathbb{Z}}{1332\mathbb{Z}}$&$\lambda(2)=2,\;\lambda(3)=\lambda(37)=1$\\

\hline
$13$&$y^2=x^3+2$&$\pm\sqrt{-11}$&$\frac{\mathbb{Z}}{2\mathbb{Z}}\times\frac{\mathbb{Z}}{6\mathbb{Z}}$&
$\frac{\mathbb{Z}}{2\mathbb{Z}}\times\frac{\mathbb{Z}}{666\mathbb{Z}}$&$\lambda(2)=2,\;\lambda(3)=\lambda(37)=1$\\

\hline
$14$&$y^2=x^3+4x-4$&$\frac{-1\pm\sqrt{-43}}{2}$&$\frac{\mathbb{Z}}{13\mathbb{Z}}$&
$\frac{\mathbb{Z}}{1343\mathbb{Z}}$&$\lambda(17)=\lambda(19)=1$\\
 \hline

\hline
$15$&$y^2=x^3+8x-4$&$-1\pm\sqrt{-10}$&$\frac{\mathbb{Z}}{14\mathbb{Z}}$&
$\frac{\mathbb{Z}}{1354\mathbb{Z}}$&$\lambda(2)=2,\;\lambda(677)=1$\\

\hline
$16$&$y^2=x^3+x+1$&$-1\pm\sqrt{-10}$&$\frac{\mathbb{Z}}{14\mathbb{Z}}$&
$\frac{\mathbb{Z}}{1354\mathbb{Z}}$&$\lambda(2)=2,\;\lambda(677)=1$\\

\hline {}&{}&{}&{}&{}&
\multicolumn{1}{|c|}{$\lambda(3)=\lambda(13)=1$}\\ \cline{6-6}
\raisebox{1.6ex}[0pt]{$17$}& \raisebox{1.6ex}[0pt]{$y^2=x^3+8x+5$}
&\raisebox{1.6ex}[0pt]{$\frac{-3\pm\sqrt{-35}}{2}$}
&\raisebox{1.6ex}[0pt]{$\frac{\mathbb{Z}}{15\mathbb{Z}}$}&
\raisebox{1.6ex}[0pt]{$\frac{\mathbb{Z}}{1365\mathbb{Z}}$}&
\multicolumn{1}{|c|}{$\lambda(5)=\lambda(7)=2$}\\

\hline {}&{}&{}&{}&{}&
\multicolumn{1}{|c|}{$\lambda(3)=\lambda(13)=1$}\\ \cline{6-6}
\raisebox{1.6ex}[0pt]{$18$}& \raisebox{1.6ex}[0pt]{$y^2=x^3+4x-1$}
&\raisebox{1.6ex}[0pt]{$\frac{-3\pm\sqrt{-35}}{2}$}
&\raisebox{1.6ex}[0pt]{$\frac{\mathbb{Z}}{15\mathbb{Z}}$}&
\raisebox{1.6ex}[0pt]{$\frac{\mathbb{Z}}{1365\mathbb{Z}}$}&
\multicolumn{1}{|c|}{$\lambda(5)=\lambda(7)=2$}\\

\hline
$19$&$y^2=x^3+8x+8$&$-2\pm\sqrt{-7}$&$\frac{\mathbb{Z}}{16\mathbb{Z}}$&
$\frac{\mathbb{Z}}{1376\mathbb{Z}}$&$\lambda(2)=2,\;\lambda(43)=1$\\

\hline $20$&$y^2=x^3+x+2$&
$-2\pm\sqrt{-7}$&$\frac{\mathbb{Z}}{2\mathbb{Z}}\times\frac{\mathbb{Z}}{8\mathbb{Z}}$&
$\frac{\mathbb{Z}}{2\mathbb{Z}}\times\frac{\mathbb{Z}}{688\mathbb{Z}}$&$\lambda(2)=2,\;\lambda(43)=1$\\

\hline

$21$&$y^2=x^3+2x+4$&$\frac{-5\pm\sqrt{-19}}{2}$&$\frac{\mathbb{Z}}{17\mathbb{Z}}$&$\frac{\mathbb{Z}}{1387\mathbb{Z}}$&$\lambda(19)=2,\;\lambda(73)=1$\\

\hline
$22$&$y^2=x^3+x+3$&$-3\pm\sqrt{-2}$&$\frac{\mathbb{Z}}{18\mathbb{Z}}$&$\frac{\mathbb{Z}}{1398\mathbb{Z}}$&$\lambda(2)=2,\;\lambda(3)=\lambda(233)=1$\\
\hline \multicolumn{6}{|c|}{\text{The\; elliptic\; curves\; of\;
numbers}\; $n$\; \text{and}\; $23-n$\; \text{are\;
twisted.}}\\\hline
\end{tabular}
\end{lrbox}
\resizebox{5.5 in}{3.5 in}{\usebox{\tablebox}}

\begin{lrbox}{\tablebox}
\begin{tabular}{|c|c|c|c|c|c|}
\multicolumn{6}{c}{\bf TABLE\;\;\; V}\\
\multicolumn{6}{c}{(\text{of\; all\; elliptic\; curves}\; E
\;\text{defined\; over}\;
$\mathbb{F}=\mathbb{Z}/13\mathbb{Z}$)}\\

\multicolumn{6}{c}{}\\
\hline

 No. &\multicolumn{1}{|c|}{$E$} & \multicolumn{1}{|c|}{\text{inverse\;roots}}&
\multicolumn{1}{|c|}{$E(\mathbb{F})$}&
\multicolumn{1}{|c|}{$K_2(E)$}&
\multicolumn{1}{|c|}{$\lambda(l)$}\\
\hline
$1$&$y^2=x^3+6$&$\frac{7\pm\sqrt{-3}}{2}$&$\frac{\mathbb{Z}}{7\mathbb{Z}}$&
$\frac{\mathbb{Z}}{2107\mathbb{Z}}$&$\lambda(7)=\lambda(43)=1$\\

\hline
$2$&$y^2=x^3+4x$&$3\pm2\sqrt{-1}$&$\frac{\mathbb{Z}}{2\mathbb{Z}}\times\frac{\mathbb{Z}}{4\mathbb{Z}}$&
$\frac{\mathbb{Z}}{2\mathbb{Z}}\times\frac{\mathbb{Z}}{1060\mathbb{Z}}$&
$\lambda(2)=2,\;\lambda(5)=\lambda(53)=1$\\
\hline

$3$&$y^2=x^3+6x+12$&$3\pm2\sqrt{-1}$&$\frac{\mathbb{Z}}{8\mathbb{Z}}$&
$\frac{\mathbb{Z}}{2120\mathbb{Z}}$&$\lambda(2)=2,\;\lambda(5)=\lambda(53)=1$\\
\hline

$4$&$y^2=x^3+3$&$\frac{5\pm3\sqrt{-3}}{2}$&$\frac{\mathbb{Z}}{3\mathbb{Z}}\times\frac{\mathbb{Z}}{3\mathbb{Z}}$&
$\frac{\mathbb{Z}}{3\mathbb{Z}}\times\frac{\mathbb{Z}}{711\mathbb{Z}}$&
$\lambda(3)=2,\;\lambda(79)=1$\\

\hline
$5$&$y^2=x^3+x+5$&$\frac{5\pm3\sqrt{-3}}{2}$&$\frac{\mathbb{Z}}{9\mathbb{Z}}$&
$\frac{\mathbb{Z}}{2133\mathbb{Z}}$&$\lambda(3)=2,\;\lambda(79)=1$\\

\hline
$6$&$y^2=x^3+2x$&$2\pm3\sqrt{-1}$&$\frac{\mathbb{Z}}{10\mathbb{Z}}$&
$\frac{\mathbb{Z}}{2146\mathbb{Z}}$&$\lambda(2)=2,\;\lambda(29)=\lambda(37)=1$\\

\hline
$7$&$y^2=x^3+8x+11$&$2\pm3\sqrt{-1}$&$\frac{\mathbb{Z}}{10\mathbb{Z}}$&
$\frac{\mathbb{Z}}{2146\mathbb{Z}}$&$\lambda(2)=2,\;\lambda(29)=\lambda(37)=1$\\

\hline
$8$&$y^2=x^3+4x+8$&$2\pm3\sqrt{-1}$&$\frac{\mathbb{Z}}{10\mathbb{Z}}$&
$\frac{\mathbb{Z}}{2146\mathbb{Z}}$&$\lambda(2)=2,\;\lambda(29)=\lambda(37)=1$\\

\hline
$9$&$y^2=x^3+8x+6$&$\frac{3\pm\sqrt{-43}}{2}$&$\frac{\mathbb{Z}}{11\mathbb{Z}}$&
$\frac{\mathbb{Z}}{2159\mathbb{Z}}$&$\lambda(17)=\lambda(127)=1$\\

\hline
$10$&$y^2=x^3+1$&$1\pm2\sqrt{-3}$&$\frac{\mathbb{Z}}{2\mathbb{Z}}\times\frac{\mathbb{Z}}{6\mathbb{Z}}$&
$\frac{\mathbb{Z}}{2\mathbb{Z}}\times\frac{\mathbb{Z}}{1086\mathbb{Z}}$&
$\lambda(2)=\lambda(3)=2,\;\lambda(181)=1$\\

\hline
$11$&$y^2=x^3+2x+5$&$1\pm2\sqrt{-3}$&$\frac{\mathbb{Z}}{12\mathbb{Z}}$&
$\frac{\mathbb{Z}}{2172\mathbb{Z}}$&$\lambda(2)=\lambda(3)=2,\;\lambda(181)=1$\\

\hline
$12$&$y^2=x^3+x+2$&$1\pm2\sqrt{-3}$&$\frac{\mathbb{Z}}{12\mathbb{Z}}$&
$\frac{\mathbb{Z}}{2172\mathbb{Z}}$&$\lambda(2)=\lambda(3)=2,\;\lambda(181)=1$\\

\hline
$13$&$y^2=x^3+8x+4$&$1\pm2\sqrt{-3}$&$\frac{\mathbb{Z}}{2\mathbb{Z}}\times\frac{\mathbb{Z}}{6\mathbb{Z}}$&
$\frac{\mathbb{Z}}{2\mathbb{Z}}\times\frac{\mathbb{Z}}{1086\mathbb{Z}}$&
$\lambda(2)=\lambda(3)=2,\;\lambda(181)=1$\\

\hline
$14$&$y^2=x^3+8x+3$&$\frac{1\pm\sqrt{-51}}{2}$&$\frac{\mathbb{Z}}{13\mathbb{Z}}$&
$\frac{\mathbb{Z}}{2185\mathbb{Z}}$&$\lambda(5)=\lambda(19)=\lambda(23)=1$\\

\hline
$15$&$y^2=x^3+x+6$&$\frac{1\pm\sqrt{-51}}{2}$&$\frac{\mathbb{Z}}{13\mathbb{Z}}$&
$\frac{\mathbb{Z}}{2185\mathbb{Z}}$&$\lambda(5)=\lambda(19)=\lambda(23)=1$\\

\hline
$16$&$y^2=x^3+x+4$&$\pm\sqrt{-13}$&$\frac{\mathbb{Z}}{14\mathbb{Z}}$&
$\frac{\mathbb{Z}}{2198\mathbb{Z}}$&$\lambda(2)=2,\;\lambda(7)=\lambda(157)=1$\\

\hline
$17$&$y^2=x^3+4x+6$&$\pm\sqrt{-13}$&$\frac{\mathbb{Z}}{14\mathbb{Z}}$&
$\frac{\mathbb{Z}}{2198\mathbb{Z}}$&$\lambda(2)=2,\;\lambda(7)=\lambda(157)=1$\\

\hline
$18$&$y^2=x^3+4x+9$&$\frac{-1\pm\sqrt{-51}}{2}$&$\frac{\mathbb{Z}}{15\mathbb{Z}}$&
$\frac{\mathbb{Z}}{2211\mathbb{Z}}$&$\lambda(3)=2,\;\lambda(11)=\lambda(67)=1$\\

\hline
$19$&$y^2=x^3+2x+2$&$\frac{-1\pm\sqrt{-51}}{2}$&$\frac{\mathbb{Z}}{15\mathbb{Z}}$&
$\frac{\mathbb{Z}}{2211\mathbb{Z}}$&$\lambda(3)=2,\;\lambda(11)=\lambda(67)=1$\\

\hline
$20$&$y^2=x^3+2x+7$&$-1\pm2\sqrt{-3}$&$\frac{\mathbb{Z}}{2\mathbb{Z}}\times\frac{\mathbb{Z}}{8\mathbb{Z}}$&
$\frac{\mathbb{Z}}{2\mathbb{Z}}\times\frac{\mathbb{Z}}{1112\mathbb{Z}}$&
$\lambda(2)=2,\;\lambda(139)=1$\\

\hline
$21$&$y^2=x^3+4x+3$&$-1\pm2\sqrt{-3}$&$\frac{\mathbb{Z}}{16\mathbb{Z}}$&$
\frac{\mathbb{Z}}{2224\mathbb{Z}}$&$\lambda(2)=2,\;\lambda(139)=1$\\

\hline
$22$&$y^2=x^3+8x+1$&$-1\pm2\sqrt{-3}$&$\frac{\mathbb{Z}}{16\mathbb{Z}}$&$
\frac{\mathbb{Z}}{2224\mathbb{Z}}$&$\lambda(2)=2,\;\lambda(139)=1$\\

\hline
$23$&$y^2=x^3+5$&$-1\pm2\sqrt{-3}$&$\frac{\mathbb{Z}}{4\mathbb{Z}}\times\frac{\mathbb{Z}}{4\mathbb{Z}}$&
$\frac{\mathbb{Z}}{4\mathbb{Z}}\times\frac{\mathbb{Z}}{556\mathbb{Z}}$&
$\lambda(2)=2,\;\lambda(139)=1$\\

\hline
$24$&$y^2=x^3+2x+4$&$\frac{-3\pm\sqrt{-43}}{2}$&$\frac{\mathbb{Z}}{17\mathbb{Z}}$&
$\frac{\mathbb{Z}}{2237\mathbb{Z}}$&$\lambda(2237)=1$\\

\hline

$25$&$y^2=x^3+x+1$&$-2\pm3\sqrt{-1}$&$\frac{\mathbb{Z}}{18\mathbb{Z}}$&
$\frac{\mathbb{Z}}{2250\mathbb{Z}}$&
$\lambda(2)=\lambda(3)=2,\;\lambda(5)=1$\\

\hline

$26$&$y^2=x^3+2x+3$&$-2\pm3\sqrt{-1}$&$\frac{\mathbb{Z}}{18\mathbb{Z}}$&
$\frac{\mathbb{Z}}{2250\mathbb{Z}}$&
$\lambda(2)=\lambda(3)=2,\;\lambda(5)=1$\\

\hline

$27$&$y^2=x^3+7x$&$-2\pm3\sqrt{-1}$&$\frac{\mathbb{Z}}{3\mathbb{Z}}\times\frac{\mathbb{Z}}{6\mathbb{Z}}$&
$\frac{\mathbb{Z}}{3\mathbb{Z}}\times\frac{\mathbb{Z}}{750\mathbb{Z}}$&
$\lambda(2)=\lambda(3)=2,\;\lambda(5)=1$\\
\hline
$28$&$y^2=x^3+4x+1$&$\frac{-5\pm3\sqrt{-3}}{2}$&$\frac{\mathbb{Z}}{19\mathbb{Z}}$&
$\frac{\mathbb{Z}}{2263\mathbb{Z}}$&$\lambda(31)=\lambda(73)=1$\\

\hline

$29$&$y^2=x^3+2$&$\frac{-5\pm3\sqrt{-3}}{2}$&$\frac{\mathbb{Z}}{19\mathbb{Z}}$&
$\frac{\mathbb{Z}}{2263\mathbb{Z}}$&
$\lambda(31)=\lambda(73)=1$\\

\hline
$30$&$y^2=x^3+8x+8$&$-3\pm2\sqrt{-1}$&$\frac{\mathbb{Z}}{20\mathbb{Z}}$&
$\frac{\mathbb{Z}}{2276\mathbb{Z}}$&$\lambda(2)=2,\;\lambda(569)=1$\\
\hline
$31$&$y^2=x^3+x$&$-3\pm2\sqrt{-1}$&$\frac{\mathbb{Z}}{2\mathbb{Z}}\times\frac{\mathbb{Z}}{10\mathbb{Z}}$&
$\frac{\mathbb{Z}}{2\mathbb{Z}}\times\frac{\mathbb{Z}}{1138\mathbb{Z}}$&
$\lambda(2)=2,\;\lambda(569)=1$\\

 \hline
$32$&$y^2=x^3+4$&$\frac{-7\pm\sqrt{-3}}{2}$&$\frac{\mathbb{Z}}{21\mathbb{Z}}$&
$\frac{\mathbb{Z}}{2289\mathbb{Z}}$&$\lambda(3)=2,\;\lambda(7)=\lambda(109)=1$\\
 \hline
 \multicolumn{6}{|c|}{\text{The\; elliptic\;
curves\; of\; numbers}\; $n$\; \text{and}\; $33-n$\; \text{are\;
twisted.}}\\ \hline

\end{tabular}
\end{lrbox}
\resizebox{5.5 in}{!}{\usebox{\tablebox}}

\end{document}